\documentclass[12pt,draft]{amsart}
\usepackage{amsmath,amsthm,latexsym,amscd,amsbsy,amssymb}
\chardef\bslash=`\\ 

\makeatletter
\def\verbatim{\interlinepenalty\@M \@verbatim
  \leftskip\@totalleftmargin\advance\leftskip2pc
  \frenchspacing\@vobeyspaces \@xverbatim}
\makeatother
\makeatletter
  \def\dgt@k{\dg@DX=-3 \dg@DY=2 \dg@SIZE=3} 
\makeatother

\makeatletter
  \def\dgt@kk{\dg@DX=3 \dg@DY=-1 \dg@SIZE=3}%
\makeatother

\theoremstyle{plain}
\newtheorem{thm}{Theorem}[section]

\newtheorem{pro}[thm]{Proposition}

\theoremstyle{definition}

\newcommand{\m}{{\mathbb M}{\mathbb E}
{\mathbb N}{\mathbb G}}

\newcommand{\h}{{\mathbb H}{\mathbb I}
{\mathbb L}{\mathbb B}{\mathbb C}}

\numberwithin{equation}{section}

\newcommand{\om}{{\mathbb O}{\mathbb P}{\mathbb M}{\mathbb E}
{\mathbb N}{\mathbb G}}

\newcommand{\on}{{\mathbb O}{\mathbb P}{\mathbb N}{\mathbb O}
{\mathbb B}{\mathbb E}{\mathbb L}}

\newcommand{\oh}{{\mathbb O}{\mathbb P}{\mathbb H}{\mathbb I}
{\mathbb L}{\mathbb B}{\mathbb C}}

\newcommand{\ohs}{{\mathbb O}{\mathbb P}{\mathbb H}{\mathbb I}
{\mathbb L}{\mathbb B}{\mathbb S}}

\begin{document}


\title[Menger (N\"{o}beling)
manifolds versus
Hilbert cube (space) manifolds -- a categorical
comparison]{Menger (N\"{o}beling)
manifolds versus
Hilbert cube (space) manifolds -- a categorical comparison}
\author{Alex Chigogidze}
\address{Department of Mathematics and Statistics,
University of Saskatche\-wan,
McLean Hall, 106 Wiggins Road, Saskatoon, SK, S7N 5E6,
Canada}
\email{chigogid@math.usask.ca}
\thanks{Both authors acknowledge the support
of their respective science foundations: NSERC and INTAS;
and the second author recalls with gratitude the
hospitality offered by the Department of Mathematics $\&$ Statistics
at the University of Saskatchewan.}

\author{V.~V.~Fedorchuk}
\address{Chair of General Topology and Geometry, Faculty
of Mechanics
and Mathematics, Moscow State University, Moscow, 119899, Russia}
\email{vitaly@gtopol.math.msu.su}

\keywords{Menger manifold, Hilbert cube manifold, $n$-homotopy,
Postnikov functor}
\subjclass{Primary: 57Q05; Secondary: 55P65}


\begin{abstract}{We show that the $n$-homotopy category
$\m_{n}$ of connected $(n+1)$-dimensional Menger manifolds is
isomorphic to the homotopy category $\h_{n}$ of connected Hilbert
cube manifolds whose $k$-dimensional homotopy groups are trivial
for each $k \geq n+1$.}
\end{abstract}

\maketitle
\markboth{A.~Chigogidze and V.~V.~Fedorchuk}{Menger (N\"{o}beling)
manifolds versus
Hilbert cube (space) manifolds -- a categorical comparison}

\section{Introduction}\label{S:intro}
Developments of Menger and Hilbert cube manifold theories made it
clear that at least from a certain point of view
the $n$-dimensional Menger compactum $\mu^{n}$ is an $n$-dimensional
analog of the Hilbert cube $Q$. Indeed if one examines the corresponding
characterization theorems \cite{bestv}, \cite{tor1} it becomes obvious
that these two spaces $\mu^{n}$ and $Q$ can not be distinguished
by means of $n$-dimensional ``tests". Moreover several authors have
noticed and emphasized strong similarities not only between the above
mentioned spaces $\mu^{n}$ and $Q$ but even between
the theories of Hilbert cube and Menger manifolds themselves
(see \cite{book} for comprehensive discussion of this topic as well as
for a complete list of references). Almost every
statement in one of these theories has a precise counterpart in the other.
Analogy extends to the non-locally compact
situation as well. Turns out that the universal $n$-dimensional
N\"{o}beling space $N_{n}^{2n+1}$ can be considered as the
$n$-dimensional counterpart of the Hilbert space $\ell_{2}$
(see \cite{ageev}, \cite{tor2}). Typical differences
between the infinite- and finite-dimensional theories
seem to have disappeared here.

Below we present an attempt to explain this phenomenom by proving
that for each $n \geq 0$ the $n$-homotopy category of
$\mu^{n+1}$manifolds is isomorphic to the
homotopy category of $Q$-manifolds whose $k$-dimensional
homotopy groups are trivial for each $k \geq n+1$.

\section{Preliminaries}\label{S:pre}
We assume
that the reader is familiar with the basics of Hilbert cube (space)
manifolds theory (see, for instance, \cite{chap}, \cite{bepe}) as well
as with the basics of Menger (N\"{o}beling) manifolds theory \cite{bestv},
\cite{book}.
All spaces in these notes are separable and metrizable.
Maps are assumed to be continuous. $\mu^{n}$ denotes the $n$-dimensional
universal Menger compactum, $N_{n}^{2n+1}$ denotes the $n$-dimensional
universal N\"{o}beling space, $Q$ denotes the Hilbert cube and the
separable Hilbert space is denoted by $\ell_{2}$. Two maps
$f,g \colon X \to Y$
are said to be
$n$-homotopic $n \geq 0$ (notation: $f \stackrel{n}{\simeq} g$)
if for any map $h \colon Z \to X$,
defined on an at most $n$-dimensional space $Z$,
the compositions $f\circ h$ and $g \circ h$ are homotopic
in the usual sense. The reader can find basic properties of
$n$-homotopic maps
in \cite{whi}, \cite{baues}, \cite{book}, \cite{chi2}. 

By $\h$ we denote the homotopy category of Hilbert
cube manifolds. For each integer $n \geq 0$ let $\h_{n}$
denote the subcategory of $\h$ objects of which have
trivial $k$-dimensional homotopy groups for each $k \geq n+1$.
Let also $\m_{n}$ stand for the $n$-homotopy category of connected
$(n+1)$-dimensional Menger manifolds.

Recall that {\em the $n$-th Postnikov Functor}
\[ {\mathcal P}_{n} \colon \text{$n$-}
\operatorname{HOMOT}(CW^{n+1}) \to
\operatorname{HOMOT}(CW)\]
 is defined as follows
(see, for instance, \cite[\S 4]{baues}). For an at most
$(n+1)$-dimensional connected $CW$-complex $X$ we obtain a $CW$-complex
${\mathcal P}_{n}(X)$ by killing $k$-dimen\-sional homotopy
groups for each $k \geq n+1$, i.e. we choose a $CW$-complex
${\mathcal P}_{n}(X)$ such that
$\left( {\mathcal P}_{n}(X)\right)^{(n+1)} = X$
and $\pi_{k}\left( {\mathcal P}_{n}(X)\right) = 0$
for each $k \geq n+1$. Since
$\pi_{k}\left( {\mathcal P}_{n}(Y)\right) = 0$ for each $k \geq n+1$,
it follows that any map $f \colon X \to Y$ between at most
$(n+1)$-dimensional connected $CW$-complexes admits an extension
${\mathcal P}_{n}(f) \colon {\mathcal P}_{n}(X) \to {\mathcal P}_{n}(Y)$.
It can easily be shown that:
\begin{itemize}
\item[$(\ast )_{n}$]
${\mathcal P}_{n}(f) \simeq {\mathcal P}_{n}(g)$
if and only if $f \stackrel{n}{\simeq} g$.
\end{itemize}

\section{Main Result}\label{S:main}
In this section we show that the categories $\m_{n}$ and $\h_{n}$ are
isomorphic. Here is a scheme used in the proof of Theorem
\ref{T:main} which allows us to produce (in a canonical way)
a Hilbert cube manifold
associated to a given $(n+1)$-dimensional Menger manifold:

\begin{multline*}
 \text{$\mu^{n+1}$-manifold} \xrightarrow{\text{triangulation}}
\text{$(n+1)$-polyhedron} \xrightarrow{\text{Postnikov functor}}\\
 \text{polyhedron with $\pi_{i} = 0$, $i \geq n+1$}
\xrightarrow{\text{multiplication by $Q$}}\\
\text{$Q$-manifold with $\pi_{i} = 0$, $i \geq n+1$.}
\end{multline*}

\begin{thm}\label{T:main}
The categories $\m_{n}$ and $\h_{n}$ are isomorphic. 
\end{thm}
\begin{proof}
We construct a functor
${\mathcal M}_{n} \colon \m_{n} \to \h_{n}$ and show
that it is an isomorphism.

{\em Definition of ${\mathcal M}_{n}$}. 

{\em Objects}. First we define the homotopy class of a connected Hilbert
cube manifold, denoted
by ${\mathcal M}_{n}\left( [M]_{n}\right)$, associated to
the $n$-homotopy type $[M]_{n}$ of a connected
$(n+1)$-dimensional
Menger manifold $M$.

Choose an at most $(n+1)$-dimensional
connected locally finite 
polyhedron $K$ and a (proper) $n$-homotopy equivalence
$\alpha \colon M \to K$ (see \cite[Proposition 5.1.3
and comment on page 103]{bestv}, \cite[Proposition 4.1.10]{book}).
Let
${\mathcal P}_{n}(K)$ denote a countable connected $CW$-complex
obtained from $K$
by killing
$k$-dimensional homotopy groups of $K$ for each $k \geq n+1$. This means that
${\mathcal P}_{n}(K)$ is a countable connected $CW$-complex such that
${\mathcal P}_{n}(K)^{(n+1)} = K$ and
$\pi_{k}\left( {\mathcal P}_{n}(K)\right) = 0$ for each $k \geq n+1$.
By \cite[Theorem 13]{whi}, there exists a homotopy
equivalence
$\gamma \colon {\mathcal P}_{n}(K) \to {\mathcal L}_{n}(K)$,
where 
 ${\mathcal L}_{n}(K)$ is a countable connected and locally finite
polyhedron. By Edwards' theorem \cite[Theorem 44.1]{chap},
the product ${\mathcal L}_{n}(K) \times Q$ is a $Q$-manifold.
Finally let
${\mathcal M}_{n}([M]_{n})$ be the homotopy type of the Hilbert
cube manifold ${\mathcal L}_{n}(K) \times Q$. Clearly,
$\pi_{k}({\mathcal L}_{n}(K) \times Q) =
\pi_{k}( {\mathcal P}_{n}(K)) = 0$ for each $k \geq n+1$.

Let us show that the homotopy type
$\left[ {\mathcal M}_{n}([M]_{n})\right]$
is well defined, i.e. it does not depend on the choices
of a Menger manifold $M$ (within the $n$-homotopy type $[M]_{n}$ of $M$),
of a polyhedron $K$, of an $n$-homotopy equivalence
$\alpha \colon M \to K$, of $CW$-complexes
${\mathcal P}_{n}(K)$, ${\mathcal L}_{n}(K)$ and
of a homotopy equivalence $\gamma$.

Indeed let $\widetilde{M}$ be an $(n+1)$-dimensional Menger
manifold and $f \colon M \to \widetilde{M}$ be
an $n$-homotopy equivalence. Let also $\widetilde{K}$ be a connected 
at most $(n+1)$-dimensional locally finite
polyhedron and $\widetilde{\alpha} \colon \tilde{M} \to \widetilde{K}$
be an $n$-homotopy equivalence. Let
${\mathcal P}_{n}(\widetilde{K})$ denote a countable
connected $CW$-complex
obtained from $\widetilde{K}$
by killing all
$k$-homotopy groups of $\widetilde{K}$ for each $k \geq n+1$ and let
$\widetilde{\gamma} \colon {\mathcal P}_{n}(\widetilde{K})
\to {\mathcal L}_{n}(\widetilde{K})$ be a homotopy equivalence
where 
 ${\mathcal L}_{n}(\tilde{K})$ is a connected locally finite polyhedron. 
Our goal is to show that the $Q$-manifolds 
${\mathcal L}_{n}(K) \times Q$
and ${\mathcal L}_{n}(\widetilde{K}) \times Q$ are homotopy equivalent.

Let $g \colon \widetilde{M} \to M$, $\beta \colon K \to M$ and
$\widetilde{\beta} \colon \widetilde{K} \to \widetilde{M}$ be
$n$-homotopy inverses of the maps $f$, $\alpha$ and
$\widetilde{\alpha}$ respectively. Observe that

\begin{multline*}
(\widetilde{\alpha}\circ f \circ \beta ) \circ
(\alpha \circ g \circ \widetilde{\beta}) = 
(\widetilde{\alpha}\circ f) \circ (\beta  \circ \alpha )
\circ (g \circ \widetilde{\beta}) \stackrel{n}{\simeq}
(\widetilde{\alpha}\circ f) \circ \operatorname{id}_{M}
\circ (g \circ \widetilde{\beta}) =\\
 \widetilde{\alpha}\circ (f \circ g) \circ \widetilde{\beta}
\stackrel{n}{\simeq} \widetilde{\alpha}\circ
\operatorname{id}_{\widetilde{M}} \circ
\widetilde{\beta} = \widetilde{\alpha}\circ
\widetilde{\beta} \stackrel{n}{\simeq}
\operatorname{id}_{\widetilde{K}}
\end{multline*}

\noindent and

\begin{multline*}
(\alpha \circ g \circ \widetilde{\beta}) \circ
(\widetilde{\alpha}\circ f \circ \beta ) 
 = 
(\alpha \circ g) \circ (\widetilde{\beta} \circ
\widetilde{\alpha}) \circ (f \circ \beta ) \stackrel{n}{\simeq}
(\alpha \circ g) \circ \operatorname{id}_{\widetilde{M}}
\circ (f \circ \beta ) =\\
 \alpha \circ (g \circ f) \circ \beta \stackrel{n}{\simeq} \alpha\circ
\operatorname{id}_{M} \circ
\beta = \alpha\circ
\beta \stackrel{n}{\simeq}
\operatorname{id}_{K}.
\end{multline*}

This shows that the composition
$\widetilde{\alpha}\circ f \circ \beta \colon K
\to \widetilde{K}$ is an $n$-homotopy equivalence.
Consequently, ${\mathcal P}_{n}(\widetilde{\alpha}\circ
f \circ \beta ) \colon {\mathcal P}_{n}(K) \to
{\mathcal P}_{n}(\widetilde{K})$ is a homotopy equivalence.
Next denote by $\delta \colon {\mathcal L}_{n}(K) \to {\mathcal P}_{n}(K)$ and $\widetilde{\delta} \colon
{\mathcal L}_{n}(\widetilde{K}) \to
{\mathcal P}_{n}(\widetilde{K})$ the homotopy inverses of $\gamma$ and
$\widetilde{\gamma}$ respectively and observe that

\begin{multline*}
 \Bigl( \gamma \circ {\mathcal P}_{n}(\alpha \circ g \circ
\widetilde{\beta}) \circ\widetilde{\delta}\Bigr) \circ
\Bigl( \widetilde{\gamma}\circ
{\mathcal P}_{n}(\widetilde{\alpha}\circ
f \circ \beta ) \circ \delta \Bigr)  =\\
\Bigl( \gamma \circ {\mathcal P}_{n}(\alpha \circ g \circ
\widetilde{\beta})\Bigr) \circ (\widetilde{\delta} \circ
\widetilde{\gamma}) \circ \Bigl( {\mathcal P}_{n}(\widetilde{\alpha}\circ
f \circ \beta ) \circ \delta \Bigr) \simeq\\
 \Bigl( \gamma \circ {\mathcal P}_{n}(\alpha \circ g
\circ \widetilde{\beta})\Bigr) \circ
\operatorname{id}_{{\mathcal P}_{n}(\widetilde{K})}\circ
\Bigl( {\mathcal P}_{n}(\widetilde{\alpha}\circ
f \circ \beta ) \circ \delta \Bigr) =  \gamma \circ
\Bigl( {\mathcal P}_{n}(\alpha \circ g \circ \widetilde{\beta})
\circ  {\mathcal P}_{n}(\widetilde{\alpha}\circ
f \circ \beta )\Bigr) \circ \delta \simeq \\
 \gamma \circ \operatorname{id}_{{\mathcal P}_{n}(K)}
\circ \delta = \gamma \circ \delta \simeq
\operatorname{id}_{{\mathcal L}_{n}(K)}.
\end{multline*}

\noindent Similarly,

\begin{multline*}
\Bigl( \widetilde{\gamma}\circ
{\mathcal P}_{n}(\widetilde{\alpha}\circ
f \circ \beta ) \circ \delta \Bigr) \circ \Bigl( \gamma \circ
{\mathcal P}_{n}(\alpha \circ g \circ
\widetilde{\beta}) \circ\widetilde{\delta}\Bigr) =\\
\Bigl( \widetilde{\gamma}\circ
{\mathcal P}_{n}(\widetilde{\alpha}\circ
f \circ \beta )\Bigr) \circ ( \delta \circ \gamma ) \circ
\Bigl( {\mathcal P}_{n}(\alpha \circ g \circ
\widetilde{\beta}) \circ\widetilde{\delta}\Bigr) \simeq\\
\Bigl( \widetilde{\gamma}\circ
{\mathcal P}_{n}(\widetilde{\alpha}\circ
f \circ \beta )\Bigr) \circ \operatorname{id}_{{\mathcal P}_{n}(K)} \circ
\Bigl( {\mathcal P}_{n}(\alpha \circ g \circ
\widetilde{\beta}) \circ\widetilde{\delta}\Bigr) = 
\widetilde{\gamma}\circ
\Bigl( {\mathcal P}_{n}(\widetilde{\alpha}\circ
f \circ \beta ) \circ
{\mathcal P}_{n}(\alpha \circ g \circ
\widetilde{\beta})\Bigr) \circ\widetilde{\delta} \simeq\\
\widetilde{\gamma}\circ
\operatorname{id}_{{\mathcal P}_{n}(\widetilde{K})}
\circ\widetilde{\delta} = \widetilde{\gamma}\circ
\widetilde{\delta} \simeq \operatorname{id}_{{\mathcal L}_{n}(\widetilde{K})}.
\end{multline*}

This shows that the composition
\[ \widetilde{\gamma}\circ {\mathcal P}_{n}(\widetilde{\alpha}\circ
f \circ \beta ) \circ \delta \colon {\mathcal L}_{n}(K) \to
{\mathcal L}_{n}(\widetilde{K})\]

\noindent is a homotopy equivalence. Then the product

\[ \Bigl( \widetilde{\gamma}\circ {\mathcal P}_{n}(\widetilde{\alpha}\circ
f \circ \beta ) \circ \delta \Bigr) \times
\operatorname{id}_{Q} \colon {\mathcal L}_{n}(K) \times Q \to
{\mathcal L}_{n}(\widetilde{K}) \times Q\]

\noindent is also a homotopy equivalent as needed. This shows that the
homotopy type $\left[ {\mathcal M}_{n}(M)\right]$ is well defined.

{\em Morphisms}. Let now $M_{1}$ and $M_{2}$ be two connected
$(n+1)$-dimensional Menger manifolds
and $f \colon M_{1} \to M_{2}$ be a map. We need to define the
homotopy class ${\mathcal M}_{n}([f]_{n})$. As before
let $\alpha_{i} \colon M_{i} \to K_{i}$ be an $n$-homotopy equivalence where
$K_{i}$ is at most
$(n+1)$-dimensional locally compact polyhedron, $i = 1,2$. Let
$\beta_{i} \colon K_{i} \to M_{i}$ denote the $n$-homotopy inverse of
$\alpha_{i}$, $i = 1,2$. The composition
$\alpha_{2}\circ f \circ \beta_{1} \colon K_{1} \to K_{2}$ admits
an extension
${\mathcal P}_{n}(\alpha_{2}\circ f \circ \beta_{1})
\colon {\mathcal P}_{n}(K_{1})
\to {\mathcal P}_{n}(K_{2})$. Next let
$\gamma_{i} \colon {\mathcal P}_{n}(K_{i}) \to
{\mathcal L}_{n}(K_{i})$ denote a homotopy equivalence
where ${\mathcal L}_{n}(K_{i})$ stands for a connected
locally compact polyhedron, $i = 1,2$. Let also
$\delta_{i} \colon {\mathcal L}_{n}(K_{i}) \to
{\mathcal P}_{n}(K_{i})$ denote the homotopy inverse of
$\gamma_{i}$, $i = 1,2$. Next consider the composition
$\gamma_{2}\circ {\mathcal P}_{n}(\alpha_{2}\circ f \circ
\beta_{1})\circ \delta_{1} \colon {\mathcal L}_{n}(K_{1})
\to {\mathcal L}_{n}(K_{2})$ and the product
\[ \Bigl(\gamma_{2}\circ {\mathcal P}_{n}(\alpha_{2}\circ f \circ
\beta_{1})\circ \delta_{1}\Bigr) \times \operatorname{id}_{Q} \colon
{\mathcal L}_{n}(K_{1}) \times Q \to {\mathcal L}_{n}(K_{2}) \times Q .\]

We let 

\[ {\mathcal M}_{n}([f]_{n}) =
\Bigl[\Bigl(\gamma_{2}\circ {\mathcal P}_{n}(\alpha_{2}\circ f \circ
\beta_{1})\circ \delta_{1}\Bigr) \times \operatorname{id}_{Q}\Bigr]\]
Finally let 
${\mathcal M}_{n}([f]_{n})$ be the homotopy class of the product map 
\[ {\mathcal P}_{n}(f) \times
\operatorname{id}_{Q} \colon {\mathcal P}_{n}(M) \times Q \to
{\mathcal P}_{n}(N) \times Q .\]

Let us show that this definition is correct. Indeed let
$g \colon M_{1} \to M_{2}$ is a map such that $f\stackrel{n}{\simeq} g$.
Then $\alpha_{2}\circ g \circ
\beta_{1} \stackrel{n}{\simeq} \alpha_{2}\circ f \circ
\beta_{1}$. As was noted earlier extensions 
\[ {\mathcal P}_{n}(\alpha_{2}\circ g \circ
\beta_{1}) , {\mathcal P}_{n}(\alpha_{2}\circ f \circ
\beta_{1}) \colon {\mathcal P}_{n}(K_{1}) \to {\mathcal P}_{n}(K_{2}) \]

\noindent are homotopic. Consequently

\[ \gamma_{2}\circ {\mathcal P}_{n}(\alpha_{2}\circ g \circ
\beta_{1})\circ \delta_{1}   \simeq \gamma_{2}\circ
{\mathcal P}_{n}(\alpha_{2}\circ f \circ
\beta_{1})\circ \delta_{1}\]

\noindent and

\[ \Bigl(\gamma_{2}\circ {\mathcal P}_{n}(\alpha_{2}\circ g \circ
\beta_{1})\circ \delta_{1}\Bigr) \times \operatorname{id}_{Q}
\simeq \Bigl(\gamma_{2}\circ {\mathcal P}_{n}(\alpha_{2}\circ f \circ
\beta_{1})\circ \delta_{1}\Bigr) \times \operatorname{id}_{Q}\]

\noindent as required. This completes the definition of
the functor ${\mathcal M}_{n}$.

Next we show that ${\mathcal M}_{n}$ is an isomorphism.

{\em Step 1}. ${\mathcal M}_{n}|\operatorname{Ob}(\m_{n})
\colon \operatorname{Ob}(\m_{n}) \to \operatorname{Ob}(\h_{n})$
is surjective.

Indeed, let $X$ be a connected $Q$-manifold such that $\pi_{k}(X) = 0$
for each $k \geq n+1$. By the triangulation theorem for
$Q$-manifolds \cite[Theorem 37.2]{chap},
there exists a locally compact polyhedron $P$ such that $X$ is
homeomorphic to the product $P \times Q$. Since the projection
$\pi_{P} \colon X = P \times Q \to P$ is a homotopy equivalence,
it follows that $\pi_{k}(P) = 0$ for each $k \geq n+1$. Let
$K = P^{(n+1)}$. According to the resolution theorem for Menger
manifolds \cite[Theorem 5.1.8 and a comment on page 102]{bestv} there exists
an $n$-homotopy equivalence $\alpha \colon M \to K$, where $M$ is a
connected $(n+1)$-dimensional Menger manifold. It easily follows from the
construction of the functor ${\mathcal M}_{n}$ that
${\mathcal M}_{n}([M]_{n}) = [X]$.

{\em Step 2}. ${\mathcal M}_{n}|\operatorname{Ob}(\m_{n})
\colon \operatorname{Ob}(\m_{n}) \to \operatorname{Ob}(\h_{n})$
is injective.

Let $M_{1}$ and $M_{2}$ be two $(n+1)$-dimensional Menger manifolds such that
${\mathcal M}_{n}([M_{1}]_{n}) = {\mathcal M}_{n}([M_{2}]_{n})$. According
to the definition of the functor ${\mathcal M}_{n}$ this means that
there exist:
\begin{enumerate}
\item
An $n$-homotopy equivalences
$\alpha_{i} \colon M_{i} \to K_{i}$, where $K_{i}$ is at most
$(n+1)$-dimensional locally compact polyhedron, $i = 1,2$.
\item
A homotopy equivalence
$\gamma_{i} \colon {\mathcal P}_{n}(K_{i})
\to {\mathcal L}_{n}(K_{i})$, where ${\mathcal L}_{n}(K_{i})$ is a
connected locally finite polyhedron, $i = 1,2$.
\item
A homotopy equivalence 
\[ r \colon \colon {\mathcal L}_{n}(K_{1})\times Q  \to
{\mathcal L}_{n}(K_{2}) \times Q .\]
\end{enumerate}

\noindent Since the projection
$\pi_{1}^{i}\colon {\mathcal L}_{n}(K_{i})\times Q \to
{\mathcal L}_{n}(K_{i})$ is a homotopy equivalence ($i = 1,2$), it follows
that there exists a homotopy equivalence
$s \colon {\mathcal L}_{n}(K_{1}) \to {\mathcal L}_{n}(K_{2})$.
Then the composition
$h = \delta_{2}\circ s \circ \gamma_{1} \colon
{\mathcal P}_{n}(K_{1}) \to {\mathcal P}_{n}(K_{2})$, where
$\delta_{2} \colon {\mathcal L}_{n}(K_{2}) \to
{\mathcal P}_{n}(K_{2})$ is the homotopy inverse of $\gamma_{2}$,
is a homotopy equivalence. Without loss of generality we may assume
that 
$h\bigl( {\mathcal P}_{n}(K_{1})^{(i)}\bigr)
\subseteq {\mathcal P}_{n}(K_{2})^{(i)}$ for each $i = 1,2,\dots$.
In particular,
\[ h(K_{1}) = h\bigl( {\mathcal P}_{n}(K_{1})^{(n+1)}\bigr)
\subseteq {\mathcal P}_{n}(K_{2})^{(n+1)} = K_{2} .\]

This defines a map $f = h|K_{1} \colon K_{1}\to K_{2}$. It is easy
to see that
${\mathcal P}_{n}(f) \simeq h$ and therefore, ${\mathcal P}_{n}(f)$ is
a homotopy equivalence. This obviously implies that $f$ is an
$n$-homotopy equivalence. Finally the composition
$\beta_{2}\circ f\circ \alpha_{1} \colon M_{1}\to M_{2}$,
where $\beta_{2} \colon K_{2} \to M_{2}$ is the $n$-homotopy
inverse of $\alpha_{2}$, is also an $n$-homotopy equivalence.
This proves that $[M_{1}]_{n} = [M_{2}]_{n}$.

{\em Step 3}. ${\mathcal M}_{n}|\operatorname{Mor}(\m_{n}) \colon
\operatorname{Mor}(\m_{n}) \to \operatorname{Mor}(\h_{n})$ is
surjective.

Let $G \colon X_{1} \to X_{2}$ be a map between connected $Q$-manifolds 
such that $\pi_{k}(X_{i}) = 0$ for each $k \geq n+1$ and $i = 1,2$. By
the triangulation theorem for $Q$-manifolds,
$X_{i} \approx P_{i} \times Q$, where $P_{i}$ is a connected locally
compact polyhedron, $i = 1,2$.
Let $j \colon P_{1} \to P_{1} \times Q$ denote a section of the projection
$\pi_{P_{1}} \colon P_{1} \times Q \to P_{1}$. Observe that
$j$ is a homotopy equivalence.
Let $F = \pi_{P_{2}}\circ G \circ j \colon P_{1} \to P_{2}$, where
$\pi_{P_{2}} \colon P_{2} \times Q \to P_{2}$ is the projection onto
the first coordinate.
Note also that $G$ is homotopic to the product map
$F \times \operatorname{id}_{Q} \colon P_{1} \times Q = X_{1} \to X_{2} =
P_{2}\times Q$. Without loss of generality we may assume that
$F\bigl( P^{(i)}_{1}\bigr) \subseteq P^{(i)}_{2}$ for each $i \geq 1$. Next consider
$(n+1)$-dimensional
Menger manifolds $M_{1}$ and $M_{2}$ and $n$-homotopy equivalences
$\alpha_{i} \colon M_{1} \to P^{(n+1)}_{i}$, $i = 1,2$. Let
$\beta_{2} \colon P_{2}^{(n+1)} \to M_{2}$ denote the
$n$-homotopy inverse of $\alpha_{2}$ and 
$f = \beta_{2}\circ \bigl( F|P_{1}^{(n+1)}\bigr)
\circ \alpha_{1} \colon M_{1} \to M_{2}$.
A straightforward verification shows that
${\mathcal M}_{n}([f]_{n}) = [F]$.

{\em Step 4}. ${\mathcal M}_{n}|\operatorname{Mor}(\m_{n}) \colon
\operatorname{Mor}(\m_{n}) \to \operatorname{Mor}(\h_{n})$
is injective.

Let $f_{i} \colon M_{1} \to M_{2}$ be two maps between
$(n+1)$-dimensional Menger manifolds such
that
${\mathcal M}_{n}([f_{1}]_{n}) = {\mathcal M}_{n}([f_{2}]_{n})$. Choose
an $n$-homotopy equivalence
$\alpha_{i} \colon M_{i} \to K_{i}$, where $K_{i}$
is at most
$(n+1)$-dimensional locally compact polyhedron, $i =1,2$. Let
$\beta_{1} \colon K_{1} \to M_{1}$ denote the
$n$-homotopy inverse of $\alpha_{1}$.
Also let $\gamma_{i} \colon {\mathcal P}_{n}(K_{i})
\to {\mathcal L}_{n}(K_{i})$ be
a homotopy equivalence, $ i = 1,2$. Finally let
$\delta_{1} \colon {\mathcal L}_{n}(K_{1}) \to
{\mathcal P}_{n}(K_{1})$ stands for the homotopy
inverse of $\gamma_{1}$.

According
to our assumption the product maps
\[ \Bigl(\gamma_{2} \circ {\mathcal P}_{n}(\alpha_{2}
\circ f_{1}\circ \beta_{1})\circ \delta_{1}\Bigr) \times
\operatorname{id}_{Q} \colon
{\mathcal L}_{n}(K_{1}) \times Q \to {\mathcal L}_{n}(K_{2}) \times Q \]

\noindent and

\[ \Bigl(\gamma_{2} \circ {\mathcal P}_{n}(\alpha_{2}
\circ f_{2}\circ \beta_{1})\circ \delta_{1}\Bigr) \times
\operatorname{id}_{Q} \colon
{\mathcal L}_{n}(K_{1}) \times Q \to {\mathcal L}_{n}(K_{2}) \times Q \]

\noindent are homotopic. This implies that ${\mathcal P}_{n}(\alpha_{2}
\circ f_{1}\circ \beta_{1}) \simeq {\mathcal P}_{n}(\alpha_{2}
\circ f_{2}\circ \beta_{1})$. By $(\ast )_{n}$,
$\alpha_{2}\circ f_{1} \circ \beta_{1}
\stackrel{n}{\simeq}\alpha_{2}\circ f_{2}\circ \beta_{1}$.
Since $\alpha_{i}$ and $\beta_{i}$, $i = 1,2$, are $n$-homotopy
equivalences it follows that $f_{1} \stackrel{n}{\simeq} f_{2}$ as required.
\end{proof}

Theorem \ref{T:main} can be restated in a slightly different
manner. First of all let $\om_{n}$ denote the category whose
objects are $n$-homotopy types of connected open subspaces of
the $(n+1)$-dimensional Menger compactum $\mu^{n+1}$ and whose
morphisms are $n$-homotopy classes of continuous maps of these
open subspaces.  Similarly let $\oh_{n}$ denote the category
whose objects are homotopy types of connected open subspaces
of the Hilbert cube $Q$ all $k$-dimensional homotopy groups of
which are trivial for each $k \geq n+1$. Morphisms of this
category are homotopy classes of continuous maps of the above
described open subspaces of $Q$.

\begin{pro}\label{P:openm}
The categories $\om_{n}$ and $\oh_{n}$ are isomorphic.
\end{pro}
\begin{proof}
In the light of Theorem \ref{T:main} it suffices to show
that the categories
$\h_{n}$ and $\oh_{n}$ as well as the categories
$\m_{n}$ and $\m_{n}$ are isomorphic.

In order to define the functor

\[ {\mathcal A}_{n} \colon \h_{n} \to \oh_{n} \]

\noindent recall \cite[Corollary 16.3]{chap} that for
any $Q$-manifold $X$ the product $X \times [0,1)$
can be embedded into the Hilbert cube as an open subspace.
Based on this observation we let

\[ {\mathcal A}_{n}([X]) = [X\times [0,1)]\;\;\text{for each}\;\; [X] \in \operatorname{Ob}(\h_{n})\]
\noindent and
\[ {\mathcal A}_{n}([f]) = [ f \times
\operatorname{id}_{[0,1)}]\;\; \text{for each}\;\;
[f] \in \operatorname{Mor}(\h_{n}).\]

\noindent A straightforward verification shows that
the functor ${\mathcal A}_{n}$ is indeed an isomorphism.

In order to define the functor 
\[ {\mathcal B}_{n} \colon \m_{n} \to \om_{n}\] 

\noindent we use concept
of the $n$-homotopy kernel $\operatorname{Ker}_{n}(M)$ of
an $(n+1)$-dimensional Menger manifold $M$
(see, \cite{chi1}, \cite[Section 4.4.1]{book}). It is important
to note
that $\operatorname{Ker}_{n}(M)$ plays the role of the product
$M \times [0,1)$ in the category of $(n+1)$-dimensional
Manger manifolds. 

Let $M$ be a connected $(n+1)$-dimensional Menger manifold.
According to \cite[Theorem 2.1]{chi1}, \cite[Theorem 4.4.3]{book} the
$n$-homotopy kernel
$\operatorname{Ker}_{n}(M)$ of $M$
can be topologically identified with an open subspace of
the $(n+1)$-dimensional universal Menger compactum $\mu^{n+1}$.
The complement $M - \operatorname{Ker}(M)$ is a $Z$-set in $M$
and consequently $M$ and $\operatorname{Ker}_{n}(M)$ are
$n$-homotopy equivalent. This allows us to let 
\[ {\mathcal B}_{n}([M]_{n}) = [\operatorname{Ker}_{n}(M)]_{n}\;\;
\text{for each}\;\; [M]_{n} \in \operatorname{Ob}(\m_{n}) .\]

Let now $f \colon M \to N$ be a map between
$\mu^{n+1}$-manifolds. Since $N - \operatorname{Ker}_{n}(N)$ is
a $Z$-set in $N$, it follows that there is a map
$f^{\prime} \colon M \to \operatorname{Ker}_{n}(N)$ as close
to $f$ as we wish (considered as maps into $N$). In particular,
we may assume \cite[Proposition 4.1.6]{book} that
$f^{\prime}$ and $f$ are $n$-homotopic (in $N$). Considerations
similar to the argument used in the proof of
\cite[Theorem2.2]{chi1}, \cite[Theorem 4.4.7]{book}
guarantee
that the $n$-homotopy class of the restriction
$f^{\prime}|\operatorname{Ker}_{n}(M) \colon
\operatorname{Ker}_{n}(M)
\to \operatorname{Ker}_{n}(N)$ of $f^{\prime}$ is uniquely
defined by the $n$-homotopy class of $f$. This observation suffices to
complete the definition of the functor ${\mathcal B}_{n}$ by letting

\[ {\mathcal B}_{n}([f]_{n}) =
[f^{\prime}|\operatorname{Ker}_{n}(M)]_{n} \;\;
\text{for each}\;\; f \colon M \to N \in \operatorname{Mor}(\m_{n}) .\]

The fact that the functor ${\mathcal B}_{n}$ is an
isomorphisms can be extracted directly from the above definition. 
\end{proof}

Another application of Theorem \ref{T:main} deals with
the categories $\on_{n}$ and $\ohs_{n}$. Objects of the category
$\on_{n}$ are {\em topological} types of connected open subspaces
of the $(n+1)$-dimensional N\"{o}beling space $N_{n+1}^{2n+3}$;
morphisms of this category are $n$-homotopy types of continuous
maps between its objects. Similarly objects of the category $\ohs_{n}$
are {\em topological} types of connected open subspaces of the
Hilbert space $\ell_{2}$ whose $k$-dimensional homotopy groups are
trivial for each $k \geq n+1$; morphisms of this category are
homotopy types of continuous maps between its objects.

\begin{pro}\label{P:openn}
The categories $\on_{n}$ and $\ohs_{n}$ are isomorphic.
\end{pro}
\begin{proof}
By Proposition \ref{P:openm}, it suffices to show
that the categories
$\on_{n}$ and $\om_{n}$ as well as the categories $\ohs_{n}$
and $\oh_{n}$ are isomorphic.

In order to define the functor 
\[ {\mathcal C}_{n} \colon \ohs_{n} \to \oh_{n} \]

\noindent we recall that the separable Hilbert space $\ell_{2}$
can be embedded into the
Hilbert cube $Q$ in such a way that the complement
$Q - \ell_{2}$ forms a ${\mathcal Z}$-skeletoid in
$Q$ \cite{bepe}. Now for a connected open subspace
$U$ of $\ell_{2}$ consider an open subspace $\widetilde{U}$ of $Q$
such that $U = \widetilde{U} \cap \ell_{2}$. It follows from the elementary
properties of ${\mathcal Z}$-skeletoids that the inclusion
$U \hookrightarrow \widetilde{U}$ is a homotopy equivalence.
Consequently, if $V$
is another open subspace of $Q$ such that $U = V \cap \ell_{2}$, then
$\widetilde{U}$ and $V$ have the same homotopy type.
This observation allows us to let

\[ {\mathcal C}_{n}(U) = [\widetilde{U}] \;\; \text{for each}\;\;
U \in \operatorname{Ob}(\ohs_{n}) .\]

For a continuous map $f \colon U \to V$ of open
subspaces of
$\ell_{2}$ consider the composition 
$\widetilde{f} \colon i_{V} \circ f \circ j_{U}$,
where $j_{U} \colon \widetilde{U} \to U$ stands
for the homotopy inverse of the inclusion
$i_{U} \colon U \to \widetilde{U}$ and
$i_{V} \colon V \to \widetilde{V}$ denotes the inclusion
 Obviously the homotopy class $[\widetilde{f}]$ is completely
determined by the homotopy class $[f]$ and
we let
\[ {\mathcal C}_{n}([f]) = [\widetilde{f}] \;\;
\text{for each}\;\; [f] \in \operatorname{Mor}(\ohs_{n}) .\]

Note that if $S$ and $T$ are homotopy
equivalent open subspaces of the Hilbert cube $Q$, then the open subspaces
$U = S \cap \ell_{2}$ and $V = T \cap \ell_{2}$ of $\ell_{2}$ are
also homotopy equivalent and according to
\cite{hescho}, \cite[Theorem 7.3]{bepe} are even topologically equivalent.
This is the only non-trivial observation involved in the verification
of the fact that the functor ${\mathcal C}_{n}$ is an
isomorphism.

We define functor 
\[ {\mathcal D}_{n}\colon \on_{n} \to \om_{n} \]

\noindent similarly. First of all we note that the
$(n+1)$-dimensional universal
Menger compactum $\mu^{n+1}$ also admits \cite{chi2} a
${\mathcal Z}$-skeletoid $\Sigma^{n+1}$ such that its complement
$\nu^{n+1} = \mu^{n+1} - \Sigma^{n+1}$ is homeomorphic
\cite{chikaty}, \cite[Theorem 5.5.5]{book}
to the $(n+1)$-dimensional universal N\"{o}beling space
$N_{n+1}^{2n+3}$. Now we proceed as above. We let
\[ {\mathcal D}_{n}(U) = [\widetilde{U}]_{n} \;\;
\text{for each}\;\; U \in \operatorname{Ob}(\on_{n}) ,\]

\noindent where $\widetilde{U}$ stands for an open subspace of
$\mu^{n+1}$ such that
$U = \widetilde{U} \cap N_{n+1}^{2n+3}$. Also we let
\[ {\mathcal D}_{n}([f]_{n}) = [\widetilde{f}]_{n}\;\;
\text{for each}\;\; [f]_{n} \in \operatorname{Mor}(\on_{n}) ,\]

\noindent where $\tilde{f}$ has the same meaning as above. The
proof of the fact that the functor ${\mathcal D}_{n}$ is an
isomorphisms is based on \cite[Propositions 5.7.5 and 5.7.8]{book} which
state that $n$-homotopy equivalent connected
open subspaces of $N_{n+1}^{2n+3}$ are homeomorphic.
\end{proof}


\end{document}